\documentclass[letterpaper, 10pt, conference]{ieeeconf}      % Use this line for a4
\IEEEoverridecommandlockouts                              % This command is only
\usepackage{amsmath,graphicx,amsfonts,amssymb,epsfig,mathrsfs}
\usepackage{subcaption}
\usepackage{tikz}
\usepackage{color}
\usepackage{rotating}
\usepackage{algorithm}
\usepackage{algpseudocode}
\usepackage{algorithmicx}
\usepackage{cite}
\usepackage{url}
\usepackage{framed}
\usepackage{bm}
\graphicspath{{./figs/},{.},{./figures}}%,
\DeclareGraphicsExtensions{.jpg,.pdf,.png,.tikz}
\usepackage{pgfplots,tikzscale}
\usetikzlibrary{plotmarks}
\usetikzlibrary{external}
\newtheorem{theorem}{Theorem}

\newtheorem{lemma}{Lemma}
\newtheorem{proposition}{Proposition}

\newcommand{\differential}{{\rm{d}}}

\newcommand{\mbx}{\ensuremath{\mathbf{x}}}
\newcommand{\mbu}{\ensuremath{\mathbf{u}}}

\newcommand{\dd}{\ensuremath{{\rm d}}}
\newcommand{\Rset}{\ensuremath{\mathbb{R}}}

\newcommand{\tset}{\ensuremath{\left[0, T\right]}}
\newcommand{\iset}{\ensuremath{i=1,2, \ldots, N}}

\renewcommand{\ae}{\text{a.e.}}
\newcommand{\supp}{\ensuremath{{\rm supp}}} % suporte

\newcommand{\xhat}{\hat{x}}

\title{\LARGE \bf
Optimal Control of Thermostatic Loads for Planning Aggregate Consumption: Characterization of Solution and Explicit Strategies
}

\author{Fernando~A.C.C.~Fontes\textsuperscript{*},
Abhishek Halder\textsuperscript{**},
Jorge Becerril\textsuperscript{*},
P.R. Kumar\textsuperscript{$\sharp$}% <-this % stops a space
\thanks{*Fernando Fontes and Jorge Becerril are with SYSTEC--ISR and the Department of Electrical and Computer Engineering, Faculty of Engineering, University of Porto, 4200--465 Porto, Portugal
 (e-mail: {\tt\small faf@fe.up.pt, jorgebecerril@fe.up.pt}).}
\thanks{**Abhishek Halder is with the Department of Applied Mathematics, University of California, Santa Cruz, CA 95064, USA
        (e-mail: {\tt\small{ahalder@ucsc.edu}}).}%
\thanks{$\sharp$ P.R. Kumar is with the Department of Electrical and Computer Engineering, Texas A\&M University, College Station, TX 77843, USA
        (e-mail: {\tt\small{prk@tamu.edu}}).}
}

\begin{document}

\maketitle
\thispagestyle{empty}
\pagestyle{empty}

%%%%%%%%%%%%%%%%%%%%%%%%%%%%%%%%%%%%%%%%%%%%%%%%%%%%%%%%%%%%%%%%%%%%%%%%%%%%%%%%
\begin{abstract}
We consider the problem of planning the aggregate energy consumption
for a set of thermostatically controlled loads for demand response, accounting price forecast trajectory and thermal comfort constraints. We address this as a continuous-time optimal control problem,
and analytically characterize the structure of its solution in the general case.
 In the special case, when the price forecast is monotone and the loads have equal dynamics, we show that it is possible to determine the solution in an explicit form.
 Taking this fact into account, we handle the non-monotone price case
 by considering several subproblems, each corresponding to a time subinterval where the price function is monotone, and then allocating to each subinterval a fraction of the total energy budget. This way,  for each time subinterval, the problem reduces to a simple convex optimization problem with a scalar decision variable, for which a descent direction is also known. The price forecasts for the day-ahead energy market typically have no more than four monotone segments, so the resulting optimization problem can be solved  efficiently with modest computational resources.
\end{abstract}

%%%%%%%%%%%%%%%
% Keywords:
%Keyword 1     Optimal Control
%Keyword 2     Energy systems
%Keyword 3     Optimization algorithms
%%%%%%%%%%%%%%%%%%%%%%%%%%%%%%%%%%%%%%%%%%%%%%%%%%%%%%%%%%%%%%%%%%%%%%%%%%%%%%%%

\section{Introduction}
\label{sec:introduction}
Thermostatically controlled loads (TCLs), such as air conditioners (ACs), are valuable as flexible resources to elicit demand response, i.e., for actively controlling the loads to offset intermittency in the generation side (e.g., due to renewables)
\cite{callaway2011achieving,bashash2011modeling,zhang2013aggregated,halder2015control,halder2017architecture}.
Utilities or load serving entities (LSEs) can dynamically exploit the \emph{thermal inertia} of the population of TCLs to strategically plan and control the aggregate consumption in a desired manner. In this paper, we consider the optimal planning problem for an LSE wherein the objective is to plan the power consumption trajectory over a time horizon to minimize the total purchase cost of energy (e.g., from a day-ahead market) while adhering to the individual thermal comfort limits and TCL dynamics constraints, given that a forecasted price trajectory is available over the planning horizon.  Here, we restrict the planning problem to single horizon case, although one can envisage solving the same in a sliding time-window manner.

Since TCLs are subject to discrete ON-OFF controls, finding and implementing the solution resulting from the optimal control subject to state (here, temperature) inequality constraints is a non-trivial task, even for simple cases (e.g., monotone price forecast).
For example, physical TCLs have minimum switching period constraints which do not allow ``holding" the TCLs at a constant temperature value over an interval of time.
This suggests accounting the switching period constraint explicitly in control design, so that the solution structure for the optimal control trajectory may become well defined. On the other hand, the computational challenge in solving mixed integer control problems brings forth the question: is it possible to recover the discrete, non-convex optimal control from the simpler convexified (albeit numerical) optimal control solution?

In Section II, we outline the optimal control problem accounting switching constraints, and describe the convex relaxation. Sections III and IV characterize structure of the solution considering general and monotone price forecasts, respectively. These results motivate a decomposition strategy (Section V) allowing us to solve simpler subproblems over monotone price segments. This paper extends our earlier results \cite{halder2019optimal} to apply Pontryagin's Maximum Principle (PMP) for the planning problem accounting switching constraints. Section VI concludes the paper.

%%%%%%%%%%%%%%%%%%%%%%%%%%%%%%%%%%%%%%%%%%%%%%%%%%%%%%%%%%%%%%%%%%%%%%%%%%%%%%%%

\section{Optimal Planning Problem}
For specificity, hereafter we refer TCLs as ACs. We consider an optimal consumption planning problem over time $t\in[0,T]$ for $N$ ACs with respective (indoor) temperature states $\{x_{i}(t)\}_{i=1}^{N}$, thermal coefficients $\{\alpha_{i},\beta_{i}\}_{i=1}^{N}$, ON-OFF controls $\{u_{i}(t)\}_{i=1}^{N}$, and initial conditions $\{x_{i0}\}_{i=1}^{N}$. We suppose that the ACs have upper and lower thermal comfort levels $\{L_{i},U_{i}\}_{i=1}^{N}$, the ambient temperature trajectory is $\hat{x}(t)>\max_{i}U_{i}$, and a total energy budget for the LSE is $E$.%
 %Let $\mathbf{u}:=(u_{1},\hdots,u_{N})^{\top}\in \{0,1\}^{N}$.%
\footnote[3]{We use the shorthand $[N]:=\{1,2,\hdots,N\}$, the abbreviation a.e. to mean ``almost everywhere",
%\footnote[4]
{ the symbol $\supp\{\cdot\}$ to denote the support of a function, and}
%\footnote
{the notation $[t]^+ := \max\{0,t\}$, and $t_{1}\wedge t_{2} := \min\{t_{1},t_{2}\}$.}
}
 Given a price forecast $\pi(t)$, assuming Newtonian thermal dynamics for indoor temperature trajectories, and that an ON AC draws power $P$, the planning problem is to minimize the energy procurement cost, i.e.,
\begin{align}
& \mbox{Minimize\ } J(\mbu) = \int_{0}^{T}  \pi(t) P \sum_{i=1}^N u_i(t) \, \dd t, \nonumber \\
& \mbox{subject to}\nonumber \\
& \dot{x}_i(t)=-\alpha_i( x_i(t) - \hat{x}(t)) - \beta_i  u_i(t),  \mbox{a.e.\ } t \in [0,T], i\in[N], \nonumber \\
& L_i \le x_i(t) \le  U_i \quad  \mbox{for all\ } t \in [0,T],  \quad i\in[N], \nonumber \\
& \int_{0}^{T}  \sum_{i=1}^{N} u_i(t) \dd t = E, \label{eq:const1}\\
& u_{i}(t) \in \{0,1\} \quad  \mbox{a.e.\ } t \in [0,T],  i\in[N]. \label{eq:const2}
\end{align}
In order to apply standard optimal control tools to characterize the solution of this planning problem, namely necessary conditions of optimality in the form of the PMP, we consider a modification to this problem, analyze its solution, and then relate the solution of the modified problem to the solution of the original one.

The main difficulties in analyzing the planning problem in its original form are constraints (\ref{eq:const1}) and (\ref{eq:const2}).
By introducing an additional state variable $x_{N+1}$, the isoperimetric constraint (\ref{eq:const1}),  can be rewritten as
$$
  \dot{x}_{N+1}(t)=u_1(t)+ u_2(t) + \ldots + u_N(t), \quad  \mbox{a.e.\ } t \in [0,T],
$$
with end-point conditions $x_{N+1}(0)=0$, $x_{N+1}(T)=E$. The difficulty with constraint (\ref{eq:const2}) is the fact that it makes the set of possible control values non-convex and  an optimal solution to this continuous-time problem might not exist.
In fact, when the optimal solution would  be to maintain the temperature constant, e.g., along
the thermal limits $U_i$ or $L_i$, the corresponding control would have to chatter between $0$ and $1$ at infinite frequency. Such a solution would not be defined when the trajectories are assumed to be measurable functions (we would have to enlarge the space of trajectories to include the so-called Young measures \cite{You69}). Also, such solution would not be practically implementable in ACs.

In order to guarantee that the optimal solution is not a chattering solution, we relax the admissible control values set to its convex hull, allowing intermediate control values,
$$
 u_i(t) \in [0,1] \quad  \mbox{a.e.\ } t \in [0,T], \quad i\in[N].
$$
A natural question arises: if the ACs only have ON-OFF control, how do we interpret and implement a solution that has intermediate control values? To address this, we note that physical AC units have a maximum on-off switching frequency (which prevents  a hypothetical chattering solution from being implemented), or equivalently a minimum switching period.

Let $T_m$ be the minimum switching period of the AC unit and
 $\hat{u}_i\in[0,1]$ be an intermediate control value.
We define an implementable equivalent control $\tilde{u}_{i}\in\{0,1\}$ to be the periodic ON-OFF control with
duty-cycle $\lambda/T_m$. In each period, we have
$$
  \tilde{u}_{i}(t)=
  	\left\{
    			\begin{array}{ll}
    			1 & \text{for}\;t \in [0,\lambda), \\
    			0 & \text{for}\;t\in [\lambda, T_m).
    			\end{array}
  	\right.
$$
The entire periodic signal
 over a time interval of length  $K T_m$, with $K$ being some positive integer,
is given by
\begin{equation} \label{eq:util}
  \tilde{u}_{i}(t)=
  	\left\{
    			\begin{array}{ll}
    			1 & \text{for}\;t \in [jT_m,jT_m+\lambda), \\
    			0 & \text{for}\;t \in [jT_m+\lambda, (j+1) T_m),
    			\end{array}
  	\right.
\end{equation}
with $j= 0,...,K-1$.
Here, the  time $\lambda$ at which the control turns ON is  such that the trajectories
of $x_i$ resulting from applying either $\hat{u}_{i}$ or $\tilde{u}_{i}$ coincide at the end of the switching period.
%\begin{center}
% \includegraphics[width=.3\textwidth]{fig1}
%\end{center}
i.e., $\lambda$ satisfies
$$
  \int_0^{T_m} \mathrm{e}^{-\alpha_i(T_m -s ) } \beta_i \hat{u}_{i}(s) \dd s =
    \int_0^{T_m} \mathrm{e}^{-\alpha_i(T_m -s ) } \beta_i \tilde{u}_{i}(s) \dd s.
$$
Assuming $\hat{x}$ is constant in that period, we have that $\hat{u}_{i}$ is also constant and we obtain
$
   \hat{u}_{i} \int_0^{T_m} \mathrm{e}^{-\alpha_i(T_m -s ) }   \dd s =
     \int_0^{\lambda} \mathrm{e}^{-\alpha_i(T_m -s ) }   \dd s,
$
implying that $\lambda$ is given explicitly by
\begin{equation} \label{eq:lambda}
    \lambda  =  \frac{1}{\alpha_i} \log (1+(\mathrm{e}^{\alpha_i T_m} - 1)\hat{u}_{i}).
\end{equation}
The state trajectory resulting from $\tilde{u}_{i}$ will over-approximate the trajectory
resulting from $\hat{u}_{i}$. When the state $x_{i}(t)$ is at the lower limit $L_{i}$,  we should instead use a control starting with OFF segment, i.e.,
for $j= 0,...,K-1$
\begin{equation} \label{eq:util2}
  \tilde{u}_{i}(t)=
  	\left\{
    			\begin{array}{ll}
    			0 & \text{for}\;t \in [jT_m,(j+1)T_m-\lambda), \\
    			1 & \text{for}\;t \in [(j+1)T_m-\lambda, (j+1)T_m),
    			\end{array}
  	\right.
\end{equation}
Omitting the scaling factor $P$ without loss of generality,  the modified problem is to minimize the energy cost over $\mathcal{U}$, the set of measurable functions $u_i: [0,T] \mapsto [0,1]$, $i\in[N]$. We refer to the following  problem as \textbf{(P)}.
\begin{align}
& \mbox{Minimize\ } J(\mbu) = \int_{0}^{T}  \pi(t)  \sum_{i=1}^N u_i(t) \, \dd t, \label{eq:obj}\\
& \mbox{subject to}\nonumber \\
& \dot{x}_i(t)=-\alpha_i( x_i(t) - \hat{x}(t)) - \beta_i  u_i(t)) \mbox{a.e.\ } t \in [0,T], i\in[N], \nonumber \\
%& \quad\quad\quad\quad\quad\quad\quad\quad\quad\quad\quad\quad\quad\quad i=1,2, \ldots, N, \\
&  \dot{x}_{N+1}(t)= \sum_{i=1}^N u_i(t), \quad  \mbox{a.e.\ } t \in [0,T],\label{eq:isop}\\
& x_i(0)=x_{i0}, \quad   i\in[N],\\
& x_{N+1}(0)=0, \quad x_{N+1}(T)=E, \\
& u_i(t) \in [0,1], \quad \quad  \mbox{a.e.\ } t \in [0,T], \quad  i\in[N], \\
& L_i \le x_i(t) \le  U_i, \quad  \mbox{for all\ } t \in [0,T], \quad  i\in[N].
\label{eq:ocp}
\end{align}
%Here, %the minimization is to be carried out over functions $\mbu \in \mathcal{U}$, with
In the next section, we analyze and characterize the solution to this problem.

\section{Characterization of solution: general case}

We start by defining two control values $\overline{u}_i, \underline{u}_i$, given by
\begin{align}
\overline{u}_i := \frac{\alpha_i}{\beta_i}(\hat{x}-U_i), \quad \underline{u}_i := \frac{\alpha_i}{\beta_i}(\hat{x}-L_i),
\label{uoverbarunderbar}	
\end{align}
that are used in the development below. These controls lead to ``zero" dynamics when the state is on each boundary, thereby permitting it to slide along the same. Specifically, the control $\overline{u}_i$ permits the state to slide along the upper boundary $U_i$, while $\underline{u}_i$ permits to slide the state along the lower boundary $L_i$.
We note that  $\overline{u}_i$ and $\underline{u}_i$ are intermediate control values and its implementation in a TCL is done using \eqref{eq:util} or \eqref{eq:util2}, respectively, together with \eqref{eq:lambda}.

The main results here require the following assumptions.

\textit{Assumption A1:}
\begin{enumerate}
\item The initial states are admissible, i.e.,
 $$
 L_i \le x_{i0} \le U_i,  \mbox{\  for all $i\in[N]$}.
 $$

\item The total energy prescribed, $E$, can be spent respecting the limits $L_i,U_i$ for all initial states, i.e., $E \in [\overline{E},\underline{E})$,
where
$
\overline{E}:=\int_0^T \sum_{i=1}^N \max\{0, \overline{u}_i(t)\} \differential t,
$
and
$
\underline{E}:=\int_0^T \sum_{i=1}^N \min\{\underline{u}_i(t),1 \} \differential t.
$

\item When the states are on the boundary of the admissible region, there is a control that drives the states into the interior of the admissible region, i.e., the values of $\alpha_i$ and $\beta_i$ are such that for all $i$, and for all possible $\hat{x}$, the temperature can rise from $L_i$ with control $u_i=0$ and can decrease from $U_i$ with control $u_i=1$
\begin{align}
 -\alpha_i (L_i-\hat{x})>0,  \quad %\label{eq:A2i}\\
 -\alpha_i (U_i-\hat{x})-\beta_i<0.  \label{eq:A2}
\end{align}

\end{enumerate}

\textit{Assumption A2:}
 The functions $\pi(\cdot)$ and $\hat{x}(\cdot)$ are differentiable.
 The function $\pi(\cdot)$ does not take the specific form
 $
   \pi(t)=A \mathrm{e}^{\alpha_i t} + B,
 $
 for some index $i\in[N]$, and some constant values $A,B$ on any subinterval of $[0,T]$ of
 nonzero measure.

 Assumption A1 guarantees the existence of at least one admissible control-state pair satisfying the constraints.
 It imposes the requirement that power of the AC unit is capable of overcoming the losses for the range of outside temperatures considered.
  Assumption A2 is of a technical nature. If a very specific growth of the price is  allowed, some algebraic coincidences lead to singular controls which are much more difficult to analyze. Assumption A2 rules out the singular control scenario.

Assumptions A1 and A2 are imposed throughout the paper.
In addition to these assumptions, for some results it is also useful to consider the following equal dynamics hypothesis H1, which permits
us to deduce further relevant properties for homogeneous populations of ACs.

  \textit{Hypothesis H1:}
  There exist constants $\alpha, \beta, L,U$ such that
  for all $i\in[N]$,
  $\alpha_i=\alpha$, $\beta_i=\beta$, $L_i=L$, $U_i=U$.

With these assumptions,  using  results from optimal control theory (see e.g. \cite{vinter_optimal_2000}), in particular applying and analyzing necessary conditions of optimality in the form of a normal maximum principle in
\cite{fontes2015normality}, we can establish the following.

%that the optimal solution for  has the following characteristics.

\begin{theorem}\label{ThmOptimalControlGeneral}
 For problem \textbf{(P)}, each component ${u^*_i}$ of the optimal control is piecewise constant, and at each time it can assume only one of the 4 values: $0$, $1$, $\overline{u}_i$, or
 $\underline{u}_i$. The value $\overline{u}_i$ occurs only when the corresponding component of the state trajectory is on the upper boundary, i.e., ${x^*_i}=U_i$, and the
  value $\underline{u}_i$ occurs only when the corresponding component of the state trajectory is on the lower boundary, i.e., ${x^*_i}=L_i$.
  Moreover,  if H1 holds then
   the transitions to the values 0 or 1 occur simultaneously for all components of the control.
\end{theorem}

\subsection{Proof of the characterization result}

First, we guarantee the existence of an optimal solution. Then we guarantee that the PMP can be written in normal form.

For $i\in[N]$, let $v_i := -\alpha_i( x_i - \hat{x}(t)) - \beta_i  u_i$, and $v_{N+1}:=\sum_{i=1}^{N}u_{i}$, where $u_i \in [0,1]$. We note that for each $(t,x)$, the set $\{(v, c) \in \Rset^{N+1} \times \Rset : c \ge \pi(t) v_{N+1}\}$, being Cartesian product of convex sets, is convex. %This can be seen by noting that this set is simply a Cartesian product of convex sets, and is therefore convex.
Combining this with assumptions A1, then problem \textbf{(P)} satisfies the conditions for existence of an optimal solution; see \cite[Thm. 23.11]{clarke_functional_2013}.
%**********************************************

That the PMP is satisfied in normal form, can be checked by  verifying that certain inward--pointing conditions are satisfied along the trajectory when the state constraint is active (see \cite{fontes_normal_2013, fontes2015normality}).
In this case, the inequalities (\ref{eq:A2}) directly imply the inward-pointing constraint qualifications guaranteeing normality.
Therefore, we can apply a strengthened version of PMP \cite[Thm. 3.2]{fontes2015normality}  and obtain the following conditions involving a scalar $\pi^*$ that can be interpreted as an intermediate price.

\begin{proposition}(\!\!\cite[Appendix A]{2019arXiv190300988F})
If $\left( \mbx^\ast, \mbu^\ast \right)$ is a local minimizer for problem \textbf{(P)},
then there exist a scalar $\pi^*$,  absolutely continuous functions $p_i,q_i: \tset \to \Rset$,  and positive Radon measures $\mu_i, \ell_i$ on $\tset$, for $i\in[N]$,  satisfying%
%\footnote[4]{We use the symbol $\supp\{\cdot\}$ to denote the support of a function.}
\begin{align}
 & \dot{p}_i(t) = \alpha_i q_i(t), \quad \ae\ t \in \tset,\\
 & q_i (t) = p_i(t) + \mu_i\{[0,t)\} -  \ell_i\{[0,t)\},  t \in [0,T),\\
  & q_i (T) = p_i(T) + \mu_i\{[0,T]\} -  \ell_i\{[0,T]\} =0,\\
& \supp \{ \ell_i \} \subset  \{ t : x_i(t) = L_i \},  \label{CS1}\\
& \supp \{ \mu_i \} \subset  \{ t : x_i(t) = U_i \},  \label{CS2}\\
& \sum_{i=1}^{N} (\pi^* -\pi(t) -  \beta_i q_i ) u_i^{*}(t) \ge \sum_{i=1}^{N} (\pi^* -\pi(t) -  \beta_i q_i ) u_i, \nonumber
\end{align}
for $\mbox{a.e.\ }t \in \tset,  u_i \in [0,1]$.
%Moreover, defining  $H[t]$ to be the Hamiltonian evaluated at time $t$,
%\begin{align*}
%  H[t]= \sum_{i=1}^N  & \alpha_i q_i  ( x^*_i(t) - \hat{x}(t)) \\
%       & + \sum_{i=1}^N (\pi^* -\pi(t) -  \beta_i q_i ) u^*_i,
%\end{align*}
%we have
%\begin{equation} \label{eq:ham-t}
%H[t]=H[0] + \int_0^t H_t[s] \dd s.
%\end{equation}

\end{proposition}

We now deduce a few lemmas, which combined together yield the result asserted in Theorem \ref{ThmOptimalControlGeneral}.

\begin{lemma}(\!\!\cite[Appendix E]{2019arXiv190300988F})
% L1
For $i\in[N]$, consider the control values $(\overline{u}_i,\underline{u}_i)$ as in (\ref{uoverbarunderbar}).
%$$= \frac{\alpha_i}{\beta_i}(\hat{x}-U_i), \mbox{\ and } \underline{u}_i= \frac{\alpha_i}{\beta_i}(\hat{x}-L_i).$$
The optimal control for problem \textbf{(P)} satisfies
  $$
  u_i^*(t)=\left\{
             \begin{array}{ll}
             1  & \mbox{\ if\ }
                \pi^* -\pi(t) -  \beta_i q_i >0, \\
              0  & \mbox{\ if\ } \pi^* -\pi(t) -  \beta_i q_i < 0,\\
             \underline{u}_i & \mbox{\ if\ } \pi^* -\pi(t) -  \beta_i q_i =0,
              x_i(t)=L_i,\\
             \overline{u}_i & \mbox{\ if\ } \pi^* -\pi(t) -  \beta_i q_i  =0,
              x_i(t)=U_i.\\
             \end{array}
           \right.
$$
\end{lemma}

It remains to analyze  whether with $\pi^* -\pi(t) -  \beta_i q_i  =0$,
other intermediate values of control could be optimal. The next lemma establishes that no intermediate control values are attained when the state is strictly within the boundaries.

 \begin{lemma}(\!\!\cite[Appendix F]{2019arXiv190300988F})
 % L2
 If Assumption 2 holds, then for any $t\in I \subset [0,T]$ such that $x_i(t) \in (L_i, U_i)$ for $i\in[N]$, the control is a piecewise constant function taking values in $\{0,1\}$.
\end{lemma}

At this point, it remains to prove the last assertion of Theorem \ref{ThmOptimalControlGeneral} concerning the synchronization of the controls when H1 holds. To this end, we proceed as follows.
 \def\tbari{\overline{t}_i}
  \def\tbarj{\overline{t}_j}
 \begin{lemma}(\!\!\cite[Appendix G]{2019arXiv190300988F})
 % L3
Assume that H1 holds. Consider two trajectories $x_i$ and $x_j$ ending in the interior of the admissible state set, i.e.,  $x_i(T),x_j(T) \in (L,U)$, with control $u_i=u_j=u_{\text{end}}$, and $u_{\text{end}}$ being either the value 0 or 1.
 Let $\tbari$ and $\tbarj$ be the respective initial instances of the maximum time interval ending in $T$ with control $u_{\text{end}}$, i.e.,
\begin{align*}
  & \tbari:=\inf \{t \in [0,T] : u_i(s)=u_{\text{end}}, s \in [t, T]\}, \\
  & \tbarj:=\inf \{t \in [0,T] : u_j(s)=u_{\text{end}}, s \in [t, T]\}.
\end{align*}
We have that $\tbari=\tbarj$.
\end{lemma}

The next lemma, whose proof uses standard dynamic programming arguments, enables us to generalize the last assertions.

 \begin{lemma}% L4
 Consider the optimal control problem \textbf{(P)} in (\ref{eq:obj})-(\ref{eq:ocp}) with solution $(\mbx^{*},\mbu^{*})$. For some given time $T^{\oplus}$ in $(0,T)$, consider the optimal control problem \textbf{(P}$^{\oplus}$\textbf{)} of minimizing
\begin{equation*}
\int_{0}^{T^{\oplus}}  \pi(t) \left(u_1(t)+ u_2(t) + \ldots + u_N(t)\right) \dd t
\end{equation*}
subject to (\ref{eq:obj})-(\ref{eq:ocp}), and
$
x_i(T^{\oplus})=x^*_i(T^{\oplus})
$
for all $i\in[N]$.
For $i\in[N]$, denote the components of the optimal state for problem \textbf{(P}$^{\oplus}$\textbf{)} as $x^{\oplus}_i(t)$. Then, $x^{\oplus}_i(t)=x^*_i(t)$ for all $t \in [0,T^{\oplus}]$, for all $i\in[N]$.
\end{lemma}

 Combining Lemma 3 and Lemma 4 by placing $T^{\oplus}$ at any instant of time for which the trajectories are in the interior of the admissible state constraint set, we conclude that all transitions of the control function to 0 or to 1, are synchronized whenever H1 holds.

\section{Characterization of the solution: monotone price case}

Consider first the case when the function $\pi$ is monotonically increasing.

\def\tilin{{t_i^{L,\text{in}}}}
\def\tiuin{{t_i^{U,\text{in}}}}
\def\tuin{{t_i^{U,\text{in}}}}
\def\tiout{t_i^{\text{out}}}
\def\tout{t^{\text{out}}}
\def\tstar{t^{*}}

\begin{proposition}(\!\!\cite[Appendix B]{2019arXiv190300988F})\label{proptstar}
 Assume that the function $\pi$ is increasing, and $\hat{x}$ is constant. Then for $i\in[N]$, there exist $\tstar_i$ such that the optimal control for problem \textbf{(P)} is
$$
   u_i^*(t)=\left\{
             \begin{array}{ll}
             1  & \mbox{\ if\ }  t < \tstar_i, \quad x_i(t) \in (L_i,U_i),\\
             \underline{u}_i & \mbox{\ if\ } t < \tstar_i, \quad x_i(t)= L_i, \\
             0  & \mbox{\ if\ } t \ge \tstar_i, \quad x_i(t) \in (L_i,U_i), \\
             \overline{u}_i & \mbox{\ if\ } t \le \tstar_i, \quad x_i(t)= U_i.
             \end{array}
           \right.
$$
Moreover, if H1 holds then for all $i,j\in[N]$ such that $i\neq j$, we have  $\tstar_i= \tstar_{j} =: t^{*}$.
\end{proposition}

%\begin{center}
% \includegraphics[angle=90,width=.3\textwidth]{fig3.pdf}
%\end{center}

In the case in which all dynamics are equal, i.e.,  H1 holds, the next result gives
the optimal solution in explicit form.

\begin{theorem}(\!\!\cite[Appendix C]{2019arXiv190300988F})
 Assume the homogeneous population hypothesis H1. Assume also that the function $\pi$ is increasing, and that $\hat{x}$ is constant.
  Then, the entry times for the temperature states at the boundary $L$, are
 \begin{equation}
 \label{tilin}
 \tilin := \frac{1}{\alpha} \log
 \frac{x_{i0} + \beta/\alpha  - \hat{x}}{L+ \beta/\alpha - \hat{x}}, \quad i\in[N];
 \end{equation}
the time needed to go from $L$ to $U$ with zero control is
 \begin{equation}
 \label{t0}
   t^0 := \frac{1}{\alpha} \log
 \frac{\hat{x}- L}{\hat{x} - U};
 \end{equation}
 and the time $\tstar$ in Proposition \ref{proptstar} solves%
% \footnote{Hereafter, we use the notation $[t]^+ := \max\{0,t\}$, and $t_{1}\wedge t_{2} := \min\{t_{1},t_{2}\}$.}
 \begin{align*}%\label{tstar}
 \!\!\displaystyle\sum_{i=1}^{N}\!\! \bigg\{\tilin\wedge\tstar  \!+\! \left[\tstar-\tilin\right]^{+}\!\!\underline{u} \!+\! \left[T - \tstar-t^0\right]^{+} \!\overline{u}\bigg\}  \!=\! E.
 \end{align*}
 Furthermore, let $\tuin:=t^0_i + \tstar$, $i\in[N]$, denote the entry times for the temperature states at the boundary $U$. In the case when $\tilin \le \tstar < \tuin \le T$ for all $i\in[N]$, the time $\tstar$  simplifies to
 $$
    \tstar = \frac{E - (1-\underline{u})\left(\sum_{i=1}^{N} \tilin\right) - N \overline{u}( T- t^0)}
                 {N(\underline{u}-\overline{u})}.
 $$
 Then, the optimal controls for problem \textbf{(P)} are
  $$
  u_i^*(t)=\left\{
             \begin{array}{ll}
             1  & \mbox{\ if\ }  t \in [0,  \tilin\wedge\tstar),\\
             \underline{u}_i & \mbox{\ if\ } t \in [\tilin\wedge\tstar, \tstar), \\
             0  & \mbox{\ if\ } t \in [\tstar, \tuin\wedge T), \\
             \overline{u}_i & \mbox{\ if\ } t \in [\tuin\wedge T, T].
             \end{array}
           \right.
$$
The corresponding optimal states are given by
 \begin{align*}
  &x_i^*(t)= \nonumber\\
  &\left\{
             \begin{array}{ll}
       e^{-\alpha t} x_{i0} + (\hat{x} -\beta/\alpha)(1 - e^{-\alpha t})
               & \mbox{\ if\ }  t \in [0, \tilin\wedge\tstar),\\
             L & \mbox{\ if\ } t \in [\tilin\wedge\tstar, \tstar), \\
             e^{-\alpha (t-\tstar)}L  & \mbox{\ if\ } t \in [\tstar,  \tuin\wedge T), \\
             U & \mbox{\ if\ } t \in [\tuin\wedge T, T].
             \end{array}
           \right.
 \end{align*}
% the  multipliers $\pi^*$ and $q$ satisfy
% \begin{align*}
% &  \dot{q}_i(t)= \alpha q_i(t) &  t \in [0, \min\{\tilin,\tstar\}) \\
% & q_i(\tilin)= \frac{1}{\beta_i}(\pi^* - \pi(\tilin)) \\
% & q_i(t)=\frac{1}{\beta_i}(\pi^* - \pi(t)) &  t \in ( \min\{\tilin,\tstar\},\tstar] \\
% &   \dot{q}_i(t)= \alpha q_i(t)                &  t \in [\tstar,\min\{\tiuin,T\}] \\
% & q_i(\tiuin)= \frac{1}{\beta_i}(\pi^* - \pi(\tiuin)) \\
% & q_i(t)= \frac{1}{\beta_i}(\pi^* - \pi(t)) &  t \in [\min\{\tiuin,T\}, T) \\
% & q_i(T)=0
% \end{align*}
%  and
%  \begin{align*}
% &  \pi^*> \pi(t)+ \beta_i q_i(t) &  t \in [0, \min\{\tilin,\tstar\})) \\
% &  \pi^*= \pi(t)+\beta_i q_i(t) & t \in [ \min\{\tilin,\tstar\}),\tstar] \\
% &  \pi^* < \pi(t)+\beta_i q_i(t) &  t \in (\tstar,\min\{\tiuin,T\}) \\
%   &  \pi^*= \pi(t)+\beta_i q_i(t) & t \in [ \min\{\tiuin,T\},T].
%  \end{align*}
If $L<x_i(T)<U$, then $\pi^*=\pi(\tstar)$; otherwise $\pi^*=\left(\pi(\tstar+t^0)-\pi(\tstar)e^{\alpha t^0}\right)/\left(1 - e^{\alpha t^0}\right)$.
\end{theorem}

The case when $\pi$ is decreasing, can be analyzed likewise. In particular,
the time needed to go from $U$ to $L$ with maximum control is
 \begin{equation}
 \label{t1}
   t^1 := \frac{1}{\alpha} \log
 \frac{U- \hat{x} + \beta/\alpha}{L- \hat{x} + \beta/\alpha}.
 \end{equation}

\subsection{Example}

Consider a planning problem for a homogeneous population with the following data:
$N=2;
T=24;
L=18;
U=22;
\xhat=30;
\alpha=0.1;
\beta=20\alpha;
P=1;
E=0.5\times N\times T=24;
\mbx_0\equiv\left(x_{10},x_{20},x_{30}\right)^{\top}=[L+1,  U-1,  0];
\pi(t)=1 +t.
$ We remind the readers that the last component of the state vector is an auxiliary state with dynamics $\dot{x}_{3} = u_{1}+u_{2}$, subject to boundary conditions $x_{3}(0)=0$, $x_{3}(T)=E$.
 From (12), we have $\bar{u}=0.4$, $\underbar{u}=0.6$.

From Theorem 2, we obtain
$\tilin = [1.1778, 3.1845]$, $t_0 =4.0547$, $\tstar =15.7469$, $\tiuin =19.8016$, the minimum value of cost $J=245.9712$, and $\pi^{*} =8.6376$.

In Fig. 1, we compare the optimal states $x_{i}^{*}(t)$, $i=1,2,3$, obtained using the numerical optimal
control solver ICLOCS \cite{falugi_iclocs_2010}, and the same obtained from Theorem 2.

%\end{document}
\begin{figure}[t]
\begin{center}
 \includegraphics[width=0.9\linewidth]{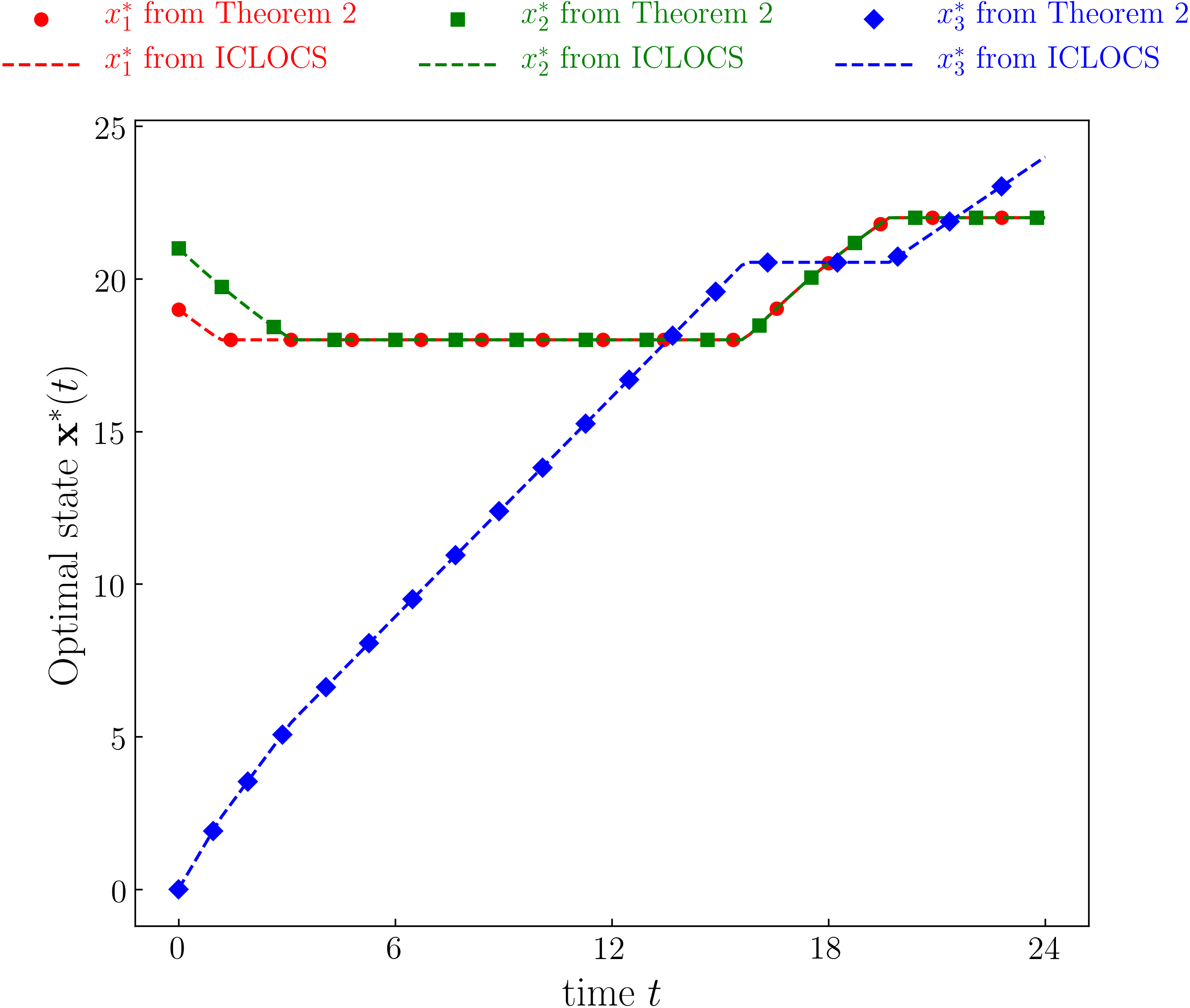}
\end{center}
\vspace*{-0.1in}
\caption{{\small{Optimal state trajectories $x_{i}^{*}(t)$, $i=1,2,3$, for the non-monotone price case with $x_{i}^{*}(0)\equiv[19, 21, 0]$, computed via the numerical optimal control solver ICLOCS, and  via the analytical solution based on Theorem 2.}}}
          \label{fig:fig1}
\vspace*{-0.25in}
\end{figure}

\section{The Nonmonotone Price Case: A Simple Explicit Strategy}

In the preceding section, we presented the solution for problem \textbf{(P)} in explicit form for the monotone price case.
In this section, we address the general price function case by decomposing \textbf{(P)} into several subproblems, each corresponding to a time subinterval where the price function is monotone.

Suppose that there are $M$ subintervals in which the price function has monotonic segments.
In practice, the price forecasts for
the day-ahead energy market typically have no more than four
monotone segments.
 Below, we describe an algorithm to iteratively compute the optimal allocation of the total energy budget in each of these $M$ subintervals. We will argue that the energy budget allocation problem can be cast as an optimization problem in $M$ scalar decision variables. The resulting problem has several features which make this approach tractable. \emph{First}, it will turn out to be a convex problem with respect to the
 $M$ scalar decision variables. \emph{Second}, the multipliers $\pi^*$, which can also be determined explicitly, define a descent direction and also an optimality criterion.
 The usefulness of the multipliers in optimal allocation of common resources has long been recognized in optimization \cite{everett_iii_generalized_1963}.

Specifically for $j=1,\hdots,M$, the multiplier $\pi_{j}^*$ in the $j$-th subinterval is precisely the multiplier associated with the corresponding isoperimetric constraint, and thus acts as a marginal cost of the energy fraction in that subinterval. The optimal solution is obtained when the $\pi_{j}^{*}$'s are all equal. So, when the values of the $\pi_{j}^*$'s are different, they help define the descent direction for optimal allocation.%, contributing to rapidly determining a solution to the problem.
They also provide a stopping criterion by detecting optimality.

\subsection{The parametric problem and its convexity}

Let the price curve have $M$ monotone segments supported over $M$ disjoint subintervals of $[0,T]$. Suppose these subintervals are of lengths $T_1$, $T_2$, \ldots $T_M$, with $T_1+T_2+ \ldots +T_M=T$. Let $\mathcal{E}:=(E_1,E_2, \ldots E_M)$ be a possible allocation of the total energy $E$ among these subintervals, i.e., $E_1+E_2+ \ldots +E_M=E$.

Consider the set of admissible allocations $\mathfrak{E}$, given by
\begin{align*}
  \mathfrak{E} := \{ & (E_1,E_2, \ldots E_M) \in  \mathbb{R}^M_+: E_1 + E_2 + \ldots E_M = E,  \\
  & E_j \in [\overline{E}_j, \underline{E}_j]
  \mbox{\ for all $m=1,2, \ldots M$} \},
\end{align*}
{ where $\overline{E}_{j} := N\overline u[T_j-t^0]^+$, and $\underline{E}_{j} := N\{(t^1\wedge T_j)+\underline u[T_j-t^1]^+\}$}
(see \cite[Appendix D]{2019arXiv190300988F} for details on these limits).
For some $\mathcal{E}=(E_1,E_2, \ldots E_M) \in \mathfrak{E}$, consider the parametric problem $ \mathcal{P}(\mathcal{E})$ given by
\begin{equation*}
\underset{\mbu\in\mathcal{U}}{\text{Minimize}} \: J(\mbu)= \!\!\int_{0}^{T}\!\!\!\!\pi(t) \left(u_1(t)+ u_2(t) + \ldots + u_N(t)\right) \dd t
%\label{eq:obj-par}
\end{equation*}
subject to
\begin{align}
& \dot{x}_i(t)=-\alpha_i( x_i(t) - \hat{x}(t)) - \beta_i  u_i(t)),  \:i\in[N], \label{eq:ocp_par1}\\
 &  \dot{x}_{N+1}(t)= \sum_{i=1}^N u_i(t), \qquad\qquad\qquad\quad  \mbox{a.e.\ } t \in [0,T],
\label{eq:isop-par}\\
& x_i(0)=x_{i0}, \qquad\qquad\qquad\qquad\qquad\;\;  i\in[N],\\
& x_{N+1}(0)=0, \\
& x_{N+1}(T_1)=E_1, \\
& x_{N+1}(T_1+T_2)=E_1+E_2, \\
& \cdots \nonumber \\
& x_{N+1}(T_1+\ldots + T_M)=E_1+\ldots + E_M,\\
& u_i(t) \in [0,1], \qquad\;\, \mbox{a.e.\ } t \in [0,T], \quad  i\in[N], \\
& L_i \le x_i(t) \le  U_i, \quad  \mbox{for all\ } t \in [0,T], \quad  i\in[N].
\label{eq:ocp_par2}
\end{align}

An important property of this problem is given in the following result.
\begin{proposition}(\!\!\cite[Appendix D]{2019arXiv190300988F})
 The parametric problem $ \mathcal{P}(\mathcal{E})$ is convex in $\mathcal{E}\in\mathfrak{E}$.
\end{proposition}
Next, we give an algorithm to solve problem \textbf{(P)} using the monotone segments of the price curve.

\subsection{Algorithm}
\begin{enumerate}
\item[(B.1)] Divide the price function $\pi(t)$ into $M$ monotone segments with the corresponding time subintervals having lengths $T_1$, $T_2$, \ldots, $T_M$, respectively.

\item[(B.2)] Choose a feasible energy allocation $\mathcal{E}\in\mathfrak{E}$, i.e., choose $E_{j} \in [\overline{E}_{j}, \underline{E}_{j}]$, $j \in [M]$, such that $\sum_{j=1}^{M}E_{j}=E$.

 \item[(B.3)] For each $j \in [M]$, compute the optimal solution as well as $\pi_{j}^{*}$ for the $j$-th segment, using the explicit formulas in Theorem 2.

 \item[(B.4)] Compare the multipliers $\pi_{j}^{*}$, $j \in [M]$, among all segments.
  \textbf{If} the multipliers $\pi_{j}^{*}$ are considered equal, \textbf{then} Stop.
   \textbf{Else} reduce $E_j$ in segments with higher $\pi^*_j$, and increase $E_j$ in segments with lower $\pi^*_j$, according to an update rule given next.

   \item[(B.5)] Repeat from (B.3).
\end{enumerate}

For iteration index $k=1,2,\hdots$, the update rule in (B.4) can be implemented as
$$
  E_{j}^{(k)}= E_{j}^{(k-1)} - \gamma \frac{ \pi^*_{j} - \widehat{\pi}^{*}}{\widehat{\pi}^{*}} {\widehat{E}}, \quad j=1,\hdots,M,
$$
where $\widehat{E}:=\sum_{j=1}^{M}E_{j}/M$, $\widehat{\pi}^{*}:=\sum_{j=1}^{M}\pi_{j}^{*}/M$, and $\gamma$ is a positive
parameter (we select $\gamma=0.5$ in the numerical example below).
{ The stopping criterion is $|\pi^*_i - \pi^*_j | < \epsilon_\pi $ for some small positive parameter $\epsilon_\pi$}.
If $E_j^{(k)} \not \in [\overline{E}_j, \underline{E}_j]$ after the update, then select the nearest extremum in this interval and rescale the remaining non-saturated $E_j$'s so that $\sum_{j=1}^{M}E_{j}=E$.
The convergence of the algorithm is sensitive to the parameter $\gamma$. If it is too small, then the convergence is slow; if it is too large, then the components of the energy allocated to each segment overshoot the average $\widehat{E}$.
 Nevertheless, being just a scalar parameter, as in the case of step-lengths in line search algorithms, it is simple to tune.

\begin{theorem}
  If  the
   multipliers resulting from the previous algorithm, at any iteration satisfies $\pi_{i}^{*}=\pi_{j}^{*}$ for all $i,j=1,\hdots,M$, then that solution is optimal.
\end{theorem}

Proof.  The concatenation of the controls, trajectories and multipliers satisfy the optimality conditions of Proposition 1. Then, the convexity of the problem guarantees optimality.

\subsection{Example}

We consider an example with non-monotone price function
$
\pi(t)=5-\sin(2 \pi t/24),
$
while keeping the remaining problem data as in Example 1.

%\begin{figure}[ht]
%\begin{center}
% \includegraphics[width=.4\textwidth]{price2} \\
%\end{center}
%\caption{Non-monotone price function (with 3 monotone segments).}
%          \label{fig:fig2}
%\end{figure}

The price function is decreasing in the segment $[0,6]$, increasing in $[6,18]$ and again decreasing in $[18,24]$.
Using the optimal control solver ICLOCS, we obtain:
cost $J=
  112.6562;
\pi^{*} =
    5.3090.
$
From the algorithm in Section V.B, we obtain that
$
E= [
    6.8054,\   12.1189,\     5.0757 ];
   J =
  112.6750.
$
The resulting optimal states $x_{i}^{*}(t)$ are illustrated in Fig. 2, both
the numerical solution obtained using the optimal control solver ICLOCS,  and  the solution obtained from the algorithm in Section V.B implementing the explicit strategy and Theorems 2-3.

We stress that the implementation of  the explicit strategy requires solving a convex optimization problem with just $M$ decision variables (in this example, the number of monotone segments is $M=3$ and the solution is almost instantaneous).
This contrasts with the situation in discrete time formulations, such as Dynamic Programming or Mixed Integer Linear Programming (MILP), where one obtains much higher dimensional problems inviting significant computational load.
For example, with a switching period $T_m$ of 1 minute, the MILP formulation has $24 \times 60 \times 2 \times N$ decision variables (see \cite[Sect. 3.3]{halder2019optimal}).

\begin{figure}[ht]
\begin{center}
 \includegraphics[width=0.85\linewidth]{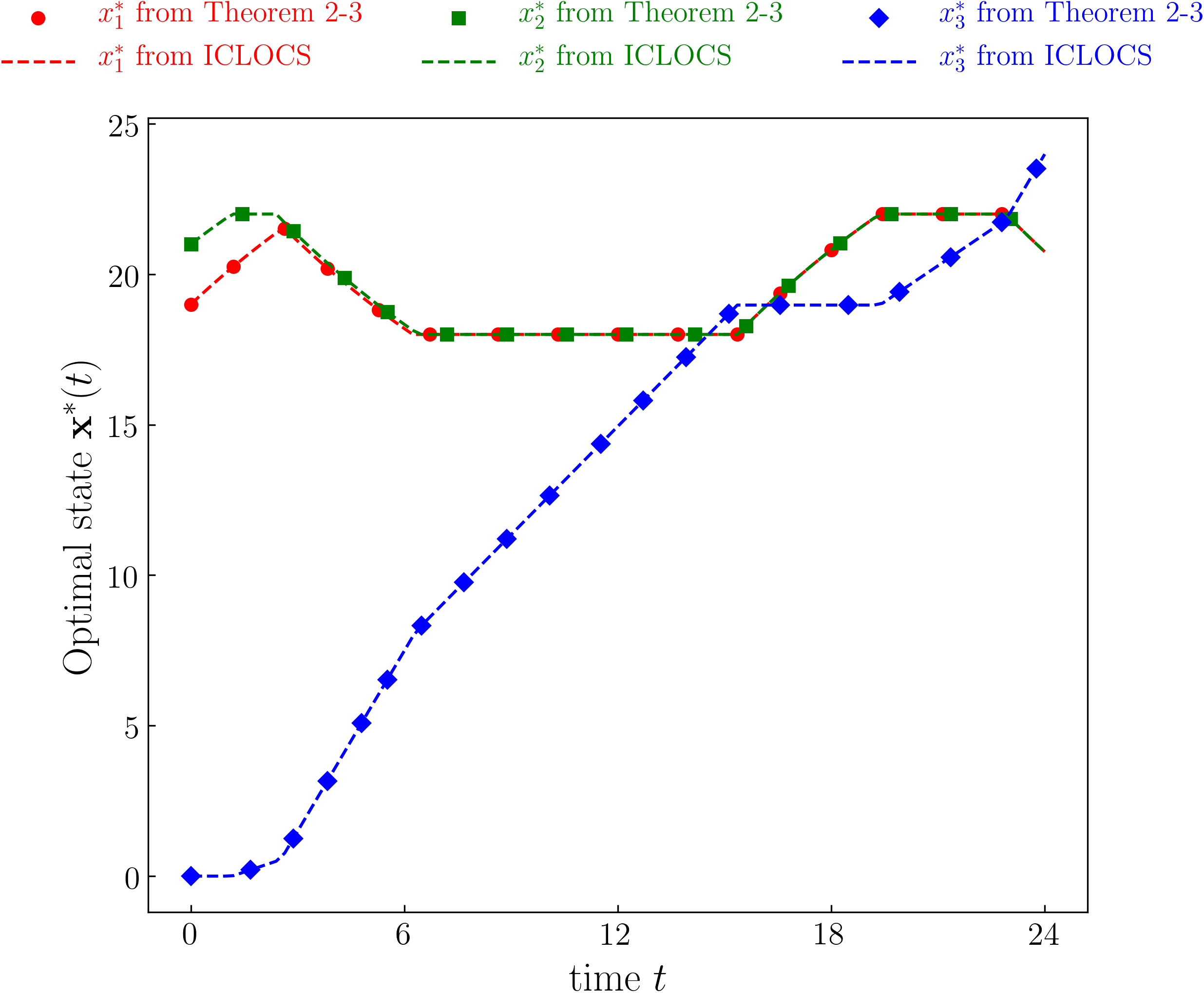}%{x2}
\end{center}
\vspace*{-0.1in}
\caption
{{\small{Optimal state trajectories, starting at $x_{i}^{*}(0)\equiv[19, 21, 0]$, computed via the numerical optimal control solver ICLOCS, and  via the proposed algorithm based on Theorems 2  and 3.}}}
          \label{fig:fig3}
\vspace*{-0.25in}
\end{figure}

\section{Conclusions}
In this paper, we considered  an optimal planning problem for demand response from the perspective of a utility or LSE, where the objective is to compute aggregate consumption for a population of thermostatic loads, conditioned on a forecasted price trajectory, that incurs the minimum cost of energy over the planning horizon. Solution of this problem can be used by the LSE for purchasing energy from the day-ahead market. A natural optimal control formulation is given that is non-convex in controls, and accounts practical switching constraints for thermostatic loads. We showed that solution of a convex relaxation can be used to recover the optimal (non-convex) solutions compliant with the switching constraints. Structural results for this relaxed problem are then exploited to further decompose this problem to sub-problems over time-intervals corresponding to monotone segments of the price forecast trajectory, which are shown to be computationally much more tractable than the original mixed-integer optimal control problem.
\section*{Acknowledgment}
{\footnotesize
This work is supported in part by NSF Science and Technology Center grant CCF-0939370, the Power Systems Engineering Research Center (PSERC), ERDF/{\-}COMPETE/{\-}NORTE2020/{\-}POCI/{\-}FCT funds through grants
  PTDC-{\-}EEI-{\-}AUT-{\-}2933-{\-}2014--Toccata,
and
02/{\-}SAICT/{\-}2017-{\-}31447-{\-}Upwind.
}
\bibliographystyle{plain}
\bibliography{OCP_bib}
%
%\end{document}
%
\newpage

\appendix
\section*{Proof of intermediate results}

The proof of some intermediate results, for lack of space, is not included in the version submitted to the IEEE Control Systems Letters and to 2019 IEEE Conference on Decision and Control. Therefore, we are placing them here for completeness and for reviewing purposes.

\subsection{Proof of Proposition 1}
%\begin{proof}
Noticing that $x_{N+1}$ is unconstrained and that (\ref{eq:isop})
does not depend on the state, the application of the maximum principle yields
$
  q_{N+1}=p_{N+1},\quad  \dot{p}_{N+1}=0.
$
That is, $q_{N+1}$ has a constant value. Denote such a value by $\pi^*$.
%
%To obtain condition (\ref{eq:ham-t}), we transform the problem into an autonomous one by considering an additional state with dynamics $\dot{t}=1$ and with initial state equal to zero. From the constancy of the Hamiltonian in the autonomous case %(\ref{eq:ham-const})
%the result (\ref{eq:ham-t}) follows.
The remaining conditions follow from direct application of Theorem 3.2 in \cite{fontes2015normality}.
%\end{proof}

\subsection{Proof of Proposition 2}

Assume the price is monotonically increasing. Then, the optimal strategy is to consume energy as early as possible while respecting the constraints. That, combined with Theorem 1 suggests the control function described as a candidate to optimal.  With the control function and the initial state defined, we can compute the trajectories and the adjoint vectors and show that the candidate solution, in fact satisfies the optimality conditions in Proposition 1.

The computations are done explicitly below, in the proof of Theorem 2, for the case H1 is satisfied.

\subsection{Proof of Theorem 2}

Assuming H1, by Thm 1 we have $\tstar_i =\tstar$ for all $i=1,2, ...,N$.
Let $\tilin$ be the entry time on the boundary $L_i$ (the first instant $t$ for which $x_i(t)=L_i$).
The first control switch might occur at $\tilin$ or $\tstar$, depending which occurs first. Therefore we have
 $$
  u_i^*(t)=\left\{
             \begin{array}{ll}
             1  & \mbox{\ if\ }  t \in [0,  \min\{\tilin,\tstar\})\\
             \underline{u}_i & \mbox{\ if\ } t \in [ \min\{\tilin,\tstar\}, \tstar).
            % 0  & \mbox{\ if\ } t \in [\tstar,  \min\{\tiuin,T\}) \\
%             \overline{u}_i & \mbox{\ if\ } t \in [\min\{\tiuin,T\}, T]
             \end{array}
           \right.
$$
Similarly, we define $\tiuin$ to be the entry time on the boundary $U_i$ (the first instant $t$ for which $x_i(t)=U_i$), and
 $$
  u_i^*(t)=\left\{
             \begin{array}{ll}
            % 1  & \mbox{\ if\ }  t \in [0,  \min\{\tilin,\tstar\})\\
%             \underline{u}_i & \mbox{\ if\ } t \in [ \min\{\tilin,\tstar\}, \tstar) \\
             0  & \mbox{\ if\ } t \in [\tstar,  \min\{\tiuin,T\}) \\
             \overline{u}_i & \mbox{\ if\ } t \in [\min\{\tiuin,T\}, T].
             \end{array}
           \right.
$$
This defines the trajectories
\begin{align*}
  &x_i^*(t)= \nonumber\\
  &\left\{
             \begin{array}{ll}
       e^{-\alpha_i t} x_{i0} + \\
       \quad (\hat{x} -\beta_i/\alpha_i)(1 - e^{-\alpha_i t})
               & \mbox{\ if\ }  t \in [0,  \min\{\tilin,\tstar\})\\
             L_i & \mbox{\ if\ } t \in [ \min\{\tilin,\tstar\}, \tstar) \\
             e^{-\alpha_i (t-\tstar)}L_i  & \mbox{\ if\ } t \in [\tstar,  \min\{\tiuin,T\}) \\
             U_i & \mbox{\ if\ } t \in [\min\{\tiuin,T\}, T],
             \end{array}
           \right.
 \end{align*}
 and, using also the optimality conditions, the adjoint multipliers satisfy
 \begin{align*}
 &  \dot{q}_i(t)= \alpha q_i(t) &  t \in [0, \min\{\tilin,\tstar\}) \\
 & q_i(\tilin)= \frac{1}{\beta_i}(\pi^* - \pi(\tilin)) \\
 & q_i(t)=\frac{1}{\beta_i}(\pi^* - \pi(t)) &  t \in ( \min\{\tilin,\tstar\},\tstar] \\
 &   \dot{q}_i(t)= \alpha q_i(t)                &  t \in [\tstar,\min\{\tiuin,T\}] \\
 & q_i(\tiuin)= \frac{1}{\beta_i}(\pi^* - \pi(\tiuin)) \\
 & q_i(t)= \frac{1}{\beta_i}(\pi^* - \pi(t)) &  t \in [\min\{\tiuin,T\}, T) \\
 & q_i(T)=0,
 \end{align*}
 where $\pi^*$ also satisfies
 \begin{align*}
 &  \pi^*> \pi(t)+ \beta_i q_i(t) &  t \in [0, \min\{\tilin,\tstar\})) \\
 &  \pi^*= \pi(t)+\beta_i q_i(t) & t \in [ \min\{\tilin,\tstar\}),\tstar] \\
 &  \pi^* < \pi(t)+\beta_i q_i(t) &  t \in (\tstar,\min\{\tiuin,T\}) \\
   &  \pi^*= \pi(t)+\beta_i q_i(t) & t \in [ \min\{\tiuin,T\},T].
  \end{align*}
  In case $x_i(T) \in (L_i,U_i)$, then $\pi^*$ is given by
   $$\pi^*=\pi(\tstar),$$
   else
   $$
     \pi^*=\frac{\pi(\tstar+t^0)-\pi(\tstar)e^{\alpha t^0}}{1 - e^{\alpha t^0}}.
   $$
   The knowledge of the trajectories enables us to compute the times $\tilin$ and $\tiuin$ explicitly.
   The time to go from $x_{i0}$ to $L_i$ with control $ u_i^*(t)=1$ is
 \begin{equation}
 \label{tilin}
 \tilin = \frac{1}{\alpha_i} \ln
 \frac{x_{i0} + \beta_i/\alpha_i  - \hat{x}}{L_i+ \beta_i/\alpha_i - \hat{x}},
 \end{equation}
 the time $t^0$ (time to go from $L_i$ to $U_i$ with zero control) is
 \begin{equation}
 \label{t0}
   t^0_i= \frac{1}{\alpha_i} \ln
 \frac{\hat{x}- L_i}{\hat{x} - U_i},
 \end{equation}
 the isoperimetric constraints impose that the time $\tstar$ solves
 \begin{align}
 \label{tstar}
   \sum_{\iset} \min\{\tilin,\tstar\}  &+ [\tstar-\tilin)]^+ \underline{u}_i \nonumber \\
   &+ [T - \tstar-t^0]^+ \overline{u}_i  = E.
 \end{align}
 In the case where $\tilin \le \tstar < \tiuin \le T$ for all $i$, $\tstar$ is given by
 the simpler expression
 $$
    \tstar = \frac{E - (1-\underline{u})\sum_{\iset} \tilin - N \overline{u}( T- t^0)}
                 {N(\underline{u}-\overline{u})},
 $$
 and
 \begin{equation}
 \label{tiuin}
   \tiuin=t^0_i + \tstar.
 \end{equation}

The case when the price function $pi$ is decreasing ia analysed in an analogous way. It involves the time needed to reach $L$ from $U$ with control $u\equiv 1$, which is of use later.
Denote by  $t^1$ the time needed to reach $L$ from $U$ with control $u\equiv 1$. Explicitly, $t^1$ is given by
\begin{align*}
t^1=\frac{1}{\alpha}\log\frac{U-\hat x+\beta/\alpha}{L-\hat x+\beta/\alpha}.
\end{align*}

\subsection{Proof of Proposition 3}
%\begin{proof}
We start by computing the minimum and maximum possible energy limits in each monotone segment $m$, respectively $\overline{E}_m$ and $\underline E_m$.

Let $0=t_0<t_1<\cdots<t_M=T$ be a partition of $[0,T]$ and $T_m:=t_m-t_{m-1}$. The minimum possible energy in the subinterval $[t_{m-1},t_m]$ is given by
\begin{align*}
\overline{E}_m=& N\left \{\int_{t_{m-1}}^{t_{m-1}+t^0}\text{d}t+\int_{t_{m-1}+t^0}^{t_m}\overline u\text{d}t\right \} \\
              =& N\overline u[t_m-t_{m-1}-t^0],
\end{align*}
when $t_{m-1}+t^0\leq t_{m}$, and is 0 when $t_{m-1}+t^0>t_{m}$. In both cases we get
\begin{align*}
\overline{E}_m=N\overline u[t_m-t_{m-1}-t^0]^+=N\overline u[T_m-t^0]^+.
\end{align*}
The maximum possible energy in the same interval is given by
\begin{align*}
\underline E_m &= N\left\{\int_{t_{m-1}}^{t_{m-1}+t^1}\text{d}t+\int_{t_{m-1}+t^1}^{t_m}\underline u\text{d}t\right\} \\
              &= N\{t^1+\underline u[t_m-t_{m-1}-t^1]\},
\end{align*}
when $t_{m-1}+t^1\leq t_m$, and
\begin{align*}
\underline E_m=N \int_{t_{m-1}}^{t_m}\text{d}t= N (t_m-t_{m-1}),
\end{align*}
when $t_{m-1}+t^1>t_m$. Thus
\begin{align*}
\underline E_m=N\{(t^1\wedge T_m)+\underline u[T_m-t^1]^+\}.
\end{align*}
In short,
\begin{align*}
\overline{E}_m&=N\overline u[T_m-t^0]^+,\\
\underline E_m&=N\{(t^1\wedge T_m)+\underline u[T_m-t^1]^+\}.
\end{align*}

Now, note that
 each state component $x_i$ can be written as an affine functional of the function $\mbu$,
since
$$
x_i(t)= \mathrm{e}^{-\alpha_i t} x_{i0} + \int_0^t
\mathrm{e}^{-\alpha_i (t-s)} (\alpha_i\hat{x}(s) - \beta_i u_i(s) ) ds.
$$
Therefore, all constraints of problem $ \mathcal{P}(\mathcal{E})$, (\ref{eq:ocp_par1}--\ref{eq:ocp_par2}), can be written in the form
\begin{align*}
  g_i(\mbu,\mathcal{E} ) = 0, \quad
  h_i(\mbu,\mathcal{E} ) \le 0,
\end{align*}
with $g_i$ and $h_i$ affine functions of $\mbu$ and $\mathcal{E}$, defining a  jointly convex domain in ($\mbu,\mathcal{E}$).
As a consequence, the set-valued map $R: \mathbb{R}^M \rightrightarrows \mathcal{U}$
$$
  R(\mathcal{E}):=\{ \mbu \in \mathcal{U}: \mbox{\ (\ref{eq:ocp_par1}--\ref{eq:ocp_par2}) are satisfied with $\mathcal{E}$} \}
$$
is convex on $\mathfrak{E}$.
That is, the set
$$
  \mathrm{Graph}(R):=\{(\mathcal{E},\mbu): \mathcal{E} \in \mathfrak{E}, \mbu \in R(\mathcal{E})  \}
$$
is convex, or equivalently \cite{fiacco_convexity_1986}, for all $\lambda \in [0,1]$, all $\mathcal{E}_1, \mathcal{E}_2 \in \mathfrak{E}$
$$
  \lambda R(\mathcal{E}_1) + (1-\lambda) R(\mathcal{E}_2) \subseteq  R(\lambda \mathcal{E}_1 + (1-\lambda)\mathcal{E}_2).
$$

 We note also that the set $\mathfrak{E}$ is convex. Therefore, using the arguments in
Prop. 2.1. in \cite{fiacco_convexity_1986}, we can show that
for all $\lambda \in [0,1]$, all $\mathcal{E}_1, \mathcal{E}_2 \in \mathfrak{E}$,
\begin{align*}
V(\lambda \mathcal{E}_1 &+ (1-\lambda) \mathcal{E}_2)
= \min_{\mbu \in R(\lambda \mathcal{E}_1+ (1-\lambda) \mathcal{E}_2)} J(\mbu) \\
&\le \min_{\mbu_1 \in R(\mathcal{E}_1),\mbu_2 \in R(\mathcal{E}_2)} J(\lambda \mbu_1 + (1-\lambda) \mbu_2) \\
&= \lambda \min_{\mbu_1 \in R(\mathcal{E}_1) } J( \mbu_1) +  (1-\lambda) \min_{\mbu_2 \in R(\mathcal{E}_2)} J(\mbu_2) \\
&=  \lambda V(\mathcal{E}_1) + (1-\lambda) V(\mathcal{E}_2).
\end{align*}
That is, $\mathcal{E} \mapsto V(\mathcal{E})$ is convex on $\mathfrak{E}$.
%\end{proof}

\subsection{Proof of Lemma 1}
%\begin{proof}
 % L1
  The maximization of the Hamiltonian condition directly yields the cases when $u_i^*(t)=0$ and when $u_i^*(t)=1$. Other intermediate values can only occur if $\pi^* -\pi(t) -  \beta_i q_i  =0$. When the trajectory is on the boundary $L_i$ for some interval of time, we must have $\dot{x}_i(t)=0$.
  The dynamic equation equal to zero  immediately yields $u_i^*(t)= \underline{u}_i= \frac{\alpha_i}{\beta_i}(\hat{x}-L_i)$. The same argument can be used on the boundary $U_i$ to get $\overline{u}_i$.
 %\end{proof}% L1
\subsection{Proof of Lemma 2}
%\begin{proof}
 % L2
  For the function $u_i$ to assume  some intermediate values not in the set $\{0,1\}$, by the maximization of the Hamiltonian we would have to have in that time interval
  $$
    \pi^* -\pi(t) -  \beta_i q_i  =0,
  $$
  %and
  $$
  \frac{d}{dt}(\pi^* -\pi(t) -  \beta_i q_i)  =0.
  $$
  Developing this last equation and substituting  $q_i$ from the previous equation, we obtain
  $$
   \frac{d}{dt} \pi(t)= \alpha_i \pi(t) - \alpha_i \pi^*.
  $$
  However, the solution of this equation is precisely of the structure that Assumption 2 rules out.
% \end{proof}% L2
\subsection{Proof of Lemma 3}
%\begin{proof}% L3
 Assume in contradiction to what we wish to prove that  $\tbari<\tbarj$,
 and, without loss of generality, that $u_{end}=1$.
 By (14)-(16), we have
 \begin{align*}
  & q_i(T)=q_j(T)=0, \\
  & \dot{q}_i(t)= \alpha q_i(t), \quad t \in [\tbarj,T],\\
  & \dot{q}_j(t)= \alpha q_j(t), \quad t \in [\tbarj,T].
 \end{align*}
 Therefore
 $$
 q_i(t)=q_j(t), \quad t \in [\tbarj,T].
 $$
 By the way $\tbarj$ is defined, we have
 \begin{align}
  & \pi^* -\pi(t) -  \beta_i q_j(t)  > 0 \quad t \in (\tbarj,T]\\
 & \pi^* -\pi(t) -  \beta_i q_j(t)  = 0 \quad t=\tbarj,
 \label{eq1}
 \end{align}
 but also
 \begin{equation}
 \pi^* -\pi(t) -  \beta_i q_i(t)  > 0 \quad t \in (\tbari,T].
 \label{eq2}
 \end{equation}
 Since $\tbari<\tbarj$, the last two equations are a contradiction.
 Repeating the same argument when $u_{end}=0$, we prove the lemma.
% \end{proof}% L3

\end{document}